\documentclass[english,12pt]{smfart}
\usepackage{amssymb}
\usepackage{amsfonts}
\usepackage{amscd}

\swapnumbers
\let\Bbb\mathbb

\def\11{{\mathbf 1}}

\def\Bbase{A_B}

\mathchardef\alphag="7C0B \mathchardef\betag="7C0C
\mathchardef\gammag="7C0D \mathchardef\deltag="7C0E
\mathchardef\varepsilong="7C22 \mathchardef\varphig="7C27
\mathchardef\psig="7C20 \mathchardef\zetag="7C10
\mathchardef\epsilong="7C0F \mathchardef\rhog="7C1A
\mathchardef\taug="7C1C \mathchardef\upsilong="7C1D
\mathchardef\iotag="7C13 \mathchardef\thetag="7C12
\mathchardef\pig="7C19 \mathchardef\sigmag="7C1B
\mathchardef\etag="7C11 \mathchardef\omegag="7C21
\mathchardef\kappag="7C14 \mathchardef\lambdag="7C15
\mathchardef\mug="7C16 \mathchardef\xig="7C18
\mathchardef\chig="7C1F \mathchardef\nug="7C17
\mathchardef\varthetag="7C23 \mathchardef\varpig="7C24
\mathchardef\varrhog="7C25 \mathchardef\varsigmag="7C26
\mathchardef\Omegag="7C0A \mathchardef\Thetag="7C02
\mathchardef\Sigmag="7C06 \mathchardef\Deltag="7C01
\mathchardef\Phig="7C08 \mathchardef\Gammag="7C00
\mathchardef\Psig="7C09 \mathchardef\Lambdag="7C03
\mathchardef\Xig="7C04 \mathchardef\Pig="7C05
\mathchardef\Upsilong="7C07

\newtheorem{theorem}[subsubsection]{Theorem}

\newtheorem{theoremm}[subsection]{Theorem}

\theoremstyle{definition}
\newtheorem{definition}[subsubsection]{Definition}

\newtheorem{def-prop}[subsubsection]{Proposition-Definition}
\newtheorem{def-theorem}[subsubsection]{Theorem-Definition}
\newtheorem{def-lem}[subsubsection]{Lemma-Definition}

\theoremstyle{remark}
\newtheorem{remark}[subsubsection]{Remark}

\theoremstyle{definition}
\newtheorem{definitionm}[subsection]{Definition}

\newtheorem{def-propm}[subsection]{Proposition-Definition}
\newtheorem{def-theoremm}[subsection]{Theorem-Definition}
\newtheorem{def-lemm}[subsection]{Lemma-Definition}

\theoremstyle{remark}

\theoremstyle{plain}

\numberwithin{equation}{subsection}

\def\boxit#1#2{\setbox1=\hbox{\kern#1{#2}\kern#1}%
\dimen1=\ht1 \advance\dimen1 by #1 \dimen2=\dp1 \advance\dimen2 by
#1
\setbox1=\hbox{\vrule height\dimen1 depth\dimen2\box1\vrule}%
\setbox1=\vbox{\hrule\box1\hrule}%
\advance\dimen1 by .4pt \ht1=\dimen1 \advance\dimen2 by .4pt
\dp1=\dimen2 \box1\relax}

\renewcommand{\theequation}{\thesubsection.\arabic{equation}}

\def\AA{{\mathbb A}}

\def\RR{{\mathbb R}}

\def\ZZ{{\mathbb Z}}

\def\cL{{\mathcal L}}
\def\cM{{\mathcal M}}

\def\cP{{\mathcal P}}

\def\11{{\mathbf 1}}

\newcommand{\comment}[1]{}
\newcommand{\Ko}{{K^\circ}}

\newcommand{\ord}{\operatorname{ord}}

\renewcommand{\epsilon}{\varepsilon}
\renewcommand{\phi}{\varphi}
\renewcommand{\emptyset}{\varnothing}

\let\Bbb\mathbb

\def\ord{{\rm ord}}

\def\Hens{{\rm Hen}}

\def\LHens{L_\Hens{ }}

\def\THens{T_\Hens{}}

\begin{document}

\title[]{An introduction to $b$-minimality}

\author{Raf Cluckers}
\address{Katholieke Universiteit Leuven,
Departement wiskunde, Celestijnenlaan 200B, B-3001 Leu\-ven,
Bel\-gium. Current address: \'Ecole Normale Sup\'erieure,
D\'epartement de ma\-th\'e\-ma\-ti\-ques et applications, 45 rue
d'Ulm, 75230 Paris Cedex 05, France} \email{cluckers@ens.fr}
\urladdr{www.wis.kuleuven.be/algebra/Raf/}

\begin{abstract}
We give a survey with some explanations but no proofs of the new
notion of $b$-minimality by the author and F. Loeser [\textit{$b$-minimality}, J. Math. Log., \textbf{7} no. 2 (2007), 195--227,
math.LO/0610183]. We compare this notion with other notions like
$o$-minimality, $C$-minimality,
$p$-minimality, and so on. 
\end{abstract}
\maketitle

\section{Introduction}

As van den Dries notes in his book \cite{vdD}, Grothendieck's dream
of tame geometries found a certain realization in model theory, at
first by the study of the geometric properties of definable sets for
some nice structure like the field of real numbers, and then by
axiomatizing these properties by notions of $o$-minimality,
minimality, $C$-minimality, $p$-minimality, $v$-minimality,
$t$-minimality, $b$-minimality, and so on. Although there is a joke
speaking of $x$-minimality with $x=a,b,c,d,\ldots,$ these notions
are useful and needed in different contexts for different kinds of
structures, for example, $o$-minimality is for ordered structures,
and $v$-minimality is for algebraically closed valued fields.

In recent work with F. Loeser \cite{CLb}, we tried to unify some of
the notions of $x$-minimality for different $x$, for certain $x$
only under extra conditions, to a very basic notion of
$b$-minimality. At the same time, we tried to keep this notion very
flexible, very tame with many nice properties, and able to describe
complicated behavior.

An observation of Grothendieck's is that instead of looking at
objects, it is often better to look at morphisms and study the
fibers of the morphisms. In one word, that is what $b$-minimality
does: while most notions of $x$-minimality focus on sets and
axiomatize subsets of the line to be simple (or tame),
$b$-minimality focuses on definable functions and gives axioms on
the existence of definable functions with nice fibers.

We give a survey on the new notion of $b$-minimality and put it in
context, without giving proofs, and refer to \cite{CLb} for the
proofs.

\section{A context}

There is a plentitude of notions of tame geometries, even just
looking at variants of $o$-minimality like quasi-$o$-minimality,
$d$-minimality, and so on. Hence there is a need for unifying
notions. Very recently, A. Wilkie \cite{Wilkie2} expanded the real
field with entire analytic functions other than $\exp$, where the
zeros are like the set of integer powers of $2$, and he shows this
structure still has a very tame geometry.
Since an $o$-minimal structure only allows finite discrete subsets
of the line, Wilkie's structure is not $o$-minimal, but it still
probably is $d$-minimal \cite{Miller}, where an expansion of the
field $\RR$ is called $d$-minimal if for every $m$ and definable
$A\subseteq\Bbb R^{m+1}$ there is some $N$ such that for all
$x\in\Bbb R^m$, the fiber $A_x:= \{y\in\Bbb R\mid (x,y)\in A\}$ of
$A$ above $x$ either has nonempty interior or is a union of $N$
discrete sets.
A similar problem exists on algebraically closed valued fields: if
one expands them with a nontrivial entire analytic function on the
line, one gets infinitely many zeros, and thus such a structure can
not be $C$-minimal nor $v$-minimal. This shows there is a need for
flexible notions of tame geometry.

\par

In \cite{CLoes}, a general theory of motivic integration is
developed, where dependence on parameters is made possible. It is
only developed for the Denef-Pas language for Henselian valued
fields (a semi-algebraic language), although only a limited number
of properties of this language are used. Hence, one needs a notion
of tame geometry for Henselian valued fields that is suitable for
motivic integration. In this paper, we give an introduction to a
notion satisfying to some extend these requirements, named
$b$-minimality, developed in detail in \cite{CLb}.

\section{$b$-minimality}\label{sbm0}

 This section is intended to sharpen the reader's intuition
before we give the formal definitions, by giving some informal
explanations on $b$-minimality. The reader who wants to see formal
definitions
first, can proceed directly with section \ref{sbm}, or go back and forth between this and the subsequent section.\\

In a $b$-minimal set-up, there are two basic kinds of sets: balls
and points. The balls are subsets of the main sort and are given by
the fibers of a single predicate $B$ in many variables, under some
coordinate projection. There are also two kinds of sorts: there is a
unique main sort and all other sorts are called auxiliary sorts
(hence there is a partitioning of the sorts in one main sort and
some auxiliary sorts). The points are just singletons. The $b$ from
$b$-minimality refers to the word balls.
\\

One sees that in any of the notions $x$-minimal with $x=o,v,C,p$, a
ball makes sense (for example, open intervals in $o$-minimal
structures), and thus $b$-minimality a priori can make sense.\\

The formal definition of $b$-minimality will be given in section
\ref{sbm}, but here we describe some reasonable and desirable
properties, of which the axioms
will be an abstraction.\\

To be in a $b$-minimal setup, a definable subset $X$ of the line $M$
in the main sort should be a disjoint union of balls and points.
Such unions might be finite, but can as well be infinite, as long as
they are ``tame" in some sense. Namely, by a tame union, we mean
that this union is the union of the fibers of a definable family,
parameterized by auxiliary parameters. Hence, infinite unions are
``allowed" as long as they are tame in this sense. To force cell
decomposition to hold, such family should be $A$-definable as soon
as $X$ is $A$-definable, with $A$ some parameters. This is the
content of the first axiom for $b$-minimality. Thus so to speak,
there is a notion of ``allowed" infinite (disjoint) union of balls
points, and any subset of the line should be such a union.
\\

Secondly, we really want balls to be different from points, and the
auxiliary sorts to be really different from the main sort. A ball
should not be a union of points (that is to say, an ``allowed" union
of points). This is captured
in the second axiom for $b$-minimality.\\

For the third axiom, the idea of a ``tame" disjoint union in a
$b$-minimal structure is needed to formulate piecewise properties.
In the third axiom, we assume a tameness property on definable
functions from the line $M$ in the main sort to $M$. Roughly, a
definable function $f:M\to M$ should be piecewise constant or
injective, where the pieces are forming a tame disjoint union, that
is, there exists a definable family whose fibers form a disjoint
union of $M$, and whose parameters are auxiliary, and on
the fibers of this family the function $f$ is constant or injective.\\

One more word on tame disjoint unions partitioning a set $X$.
Instead of speaking of ``a" definable family whose fibers form a
partition of $X$ and whose parameters are auxiliary, we will just
speak of a definable function
$$
f:X\to S
$$
with $S$ auxiliary, and the fibers of $f$ then form such a tame
union.

\section{$b$-minimality: the definition}\label{sbm}

\subsection{Some conventions} All languages will have a unique main sort, the other
sorts are auxiliary sorts. An expansion of a language may introduce
new auxiliary sorts. If a model is named $\cM$, then the main sort
of $\cM$ is denoted by $M$.

 By \emph{definable} we shall always  mean definable
with parameters, as opposed to $\cL(A)$-definable or $A$-definable,
which means definable with parameters in $A$. By a \emph{point} we
mean a singleton. A definable set is called \emph{auxiliary} if it
is a subset of a finite Cartesian product of (the universes of)
auxiliary sorts.

If $S$ is a sort, then its Cartesian power $S^0$ is considered to be
a point and to be $\emptyset$-definable.

Recall that $o$-minimality is about expansions of the language
$\cL_{<}$ with one predicate $<$, with the requirement  that the
predicate $<$ defines a dense linear order without endpoints. In the
present setting we shall  study expansions of a language $\cL_B$
consisting of one predicate $B$, which is nonempty and which has
fibers in the $M$-sort (by definition called balls). In both
instances of tame geometry, the expansion has to satisfy extra
properties. A priori, it is not determined to which product of sorts
the predicate $B$ corresponds; this will always be fixed by the
context, or it will be supposed to be fixed later on by some
context, when it needs to be fixed.

\subsection{}
Let $\cL_B$ be the language with one predicate $B$. We require that
$B$ is interpreted in any $\cL_B$-model $\cM$ with main sort $M$ as
a nonempty set $B(\cM)$ with
 $$B(\cM)\subset \Bbase \times M$$ where $\Bbase $ is a
finite Cartesian product of (the universes of) some of the sorts of
$\cM$.

When  $a\in \Bbase $ we write $B(a)$ for
$$
B(a):=\{m\in M\mid (a,m)\in B(\cM)\},
$$
and if $B(a)$ is nonempty, we call it a \emph{ball} (in the
structure $\cM$), or $B$-ball when useful.

\begin{definition}[$b$-minimality]\label{bm}
Let $\cL$ be any expansion of $\cL_B$. We call an $\cL$-model $\cM$
\emph{$b$-minimal} when the following three conditions are satisfied
for every  set of parameters $A$ (the elements of $A$ can belong to
any of the sorts), for every $A$-definable subset $X$ of $M$, and
for every $A$-definable function $F:X\to M$.
\begin{enumerate}
 \item[(b1)]
There exists a $A$-definable function $ f:X\to S$ with $S$ an auxiliary set
such that for each $s\in f(X)$ the fiber $f^{-1}(s)$ is a point or a
ball.
 \item[(b2)] If $g$ is a definable function from an auxiliary set to
 a ball, then $g$ is not surjective.
  \item[(b3)] There exists a $A$-definable function $f : X
\rightarrow S$ with $S$ an auxiliary set such that for each $s\in f(X)$ the
restriction
 $F_{|f^{-1}(s)}$ is either injective or constant.
\end{enumerate}
We call an $\cL$-theory \emph{$b$-minimal} if all its models are
$b$-minimal.
\end{definition}

\section{Cell decomposition}

In his paper on decision procedures, Cohen \cite{Cohen} develops
cell decomposition techniques for real and $p$-adic fields, by a
kind of Taylor approximation of roots of polynomials. At that time,
the writing was rather complicated and it was only through the work
by Denef \cite{Denef3}\cite{Denef2} that some concrete notion of
$p$-adic cells became apparent. One should keep in mind that there
was no ideological framework of $o$-minimality which later on formed
intuition of what cells should be and what they should do. An
example of a fracture with actual $o$-minimal intuition about cells
was that these original $p$-adic cells were not literally designed
to partition definable sets into cells, they merely helped to
partion into some nice pieces. On these nice pieces, one could get
good properties of functions defined on them, which helped to
calculate $p$-adic integrals
\cite{Denef3}\cite{Denef2}\cite{Pas}\cite{Pas2}.\\

Another aspect of $o$-minimal intuition is that cells in one
variable should be simple and defined by induction on the variables,
both aspects were not so clear for the original $p$-adic cells and
became even more complicated in the Pas-framework. Also cell
decomposition for $C$-minimal structures \cite{HM} is somehow
complicated. In $v$-minimality \cite{UdiKazh}, cell decomposition
appears mainly implicitely.\\

The notion of $b$-minimality is intended to give a blueprint for a
versatile kind of cell decomposition for tame geometries that is
simple in one variable and defined by induction on the variables. A
(1)-cell is a tame union of balls, and a (0)-cell is a tame union of
points. Then one builds further with more variables.

\section{Cell decomposition: the definitions}
Let $\cL$ be any expansion of $\cL_B$, as before, and let $\cM$ be
an $\cL$-model.

\begin{definition}[Cells]
If all fibers of some $f:X\to S$ as in (b1) are balls, then call
$(X,f)$ a $(1)$-cell with presentation $f$. If all fibers of $f$ as
in (b1) are points, then call $(X,f)$ a $(0)$-cell with presentation
$f$. For short, call such $X$ a cell.

Let $X\subset M^n$ be definable and let $(j_1,\ldots,j_n)$ be in
$\{0,1\}^n$. Let $p_n:X\to M^{n-1}$ be the projection. Call $X$ a
$(j_1,\ldots,j_n)$-cell with presentation
$$
f:X\to S
$$
for some auxiliary $S$, when for each $\hat
x:=(x_1,\ldots,x_{n-1})\in p_n(X)$, the set $p_n^{-1}(\hat x)\subset
M$ is a $(j_n)$-cell with presentation
$$
p_n^{-1}(\hat x)\to S: x_n\mapsto f(\hat x,x_n)
$$
and $p_n(X)$ is a $(j_1,\ldots,j_{n-1})$-cell with presentation
$$
f':p_n(X)\mapsto S'
$$
for some  $f'$ satisfying  $f'\circ p_n=p\circ f$ for some $p:S\to
S'$.
%
%
%
\end{definition}

One proves that if $X$ is a $(i_1,\ldots,i_n)$-cell, then $X$ is not
a $(i'_1,\ldots,i'_n)$-cell, for the same ordering of the factors of
$M^n$,  for any tuple $(i'_1,\ldots,i'_n)$ different from
$(i_1,\ldots,i_n)$. Thus $(i_1,\ldots,i_n)$ can be called the
\emph{type}
of the $(i_1,\ldots,i_n)$-cell $X$.\\

One proves the cell decomposition theorem by compactness.

\begin{theoremm}[Cell decomposition]\label{ncd}
Let $\cM$ be a model of a $b$-minimal theory. Let $X\subset M^n$ be
a definable set. Then there exists a finite partition of $X$ into
cells.
\end{theoremm}

\section{Refinements}
Often, one has a cell decomposition of $X$, but one needs a finer
cell decomposition, such that more properties hold on the parts.
Here, it is not only the cells $X_i$ that should be partitioned
further into cells, but each $X_i$ is already written as a union of
fibers which resemble products of balls and points, and all these
fibers should be partitioned into finer parts to speak of a genuine
refinement.

\begin{definition}
Let $\cP$ and $\cP'$  be two finite partitions of $X$ into cells
$(X_i,f_i)$, resp.~$(Y_j,g_j)$. Call $\cP'$ a \emph{refinement of
$\cP$} when for each $i$ there exists $j$ such that $ Y_j\subset X_i
$ and such that $g_j$ is a refinement of $f_{ij}:=f_{i|Y_j}$, that
is, for each $a\in g_j(Y_j)$, there exists a (necessarily unique)
$b\in
 f_{ij}(Y_j)$ such that
$$
  g_j^{-1}(a)\subset f_{ij}^{-1}(b).
$$
\end{definition}

One proves by compactness that refinements exist.

\section{Relative cells}

Cells use an order of the variables, so they are very well suited to
work relatively over some of the variables.

Since in a $b$-minimal set-up there are many sorts, not all
definable sets are subsets of $M^n$, with $M$ the main sort. Still,
we want most notions to make sense for the main sort, and not to
bother about the auxiliary sorts, as long as (b1), (b2) and (b3) are
not violated. So there is a need to define all the concepts for
definable subsets of $S\times M^n$ with $S$ auxiliary, or more
generally, for definable subsets $X$ of $Y\times M^n$ for any
definable $Y$. That way, one defines relative dimension over $Y$,
cells over $Y$, a presentation over $Y$, and so on.

We just define a $(i)$-cell over $Y$ with $i=0,1$. A definable set
$X\subset Y\times M$ is called a $i$-cell with presentation
$$
f:X\to Y\times S
$$
with $S$ auxiliary if $f$ commutes with the projections $Y\times
S\to Y$ and $p:X\to Y$ to $Y$, and for each $y$, the set $p^{-1}(y)$
is a $(i)$-cell with presentation
$$
p^{-1}(y)\to S:m\mapsto f(y,m),
$$
where we have identified $\{y\}\times S$ with $S$ and $p^{-1}(y)$
with a subset of $M$.

\section{Dimension theory}\label{sdim}

Very similar to the $o$-minimal dimension as in \cite{vdD}, a
dimension theory for $b$-minimal structures unfolds.

There are many sorts, but we want the dimension to live in the main
sort.

\begin{definition}
 The dimension of a nonempty
definable set $X\subset M^n$ is defined as the maximum of all sums
$$i_1+\ldots +i_n$$
where $(i_1,\ldots,i_n)$ runs over the types of all cells contained
in $X$, for all orderings of the $n$ factors of $M^n$. To the empty
set we assign the dimension $-\infty$.

If $X\subset S\times M^{n}$ is definable with $S$ auxiliary, the
dimension of $X$ is defined as the dimension of $p(X)$ with
$p:S\times M^{n}\to M^n$ the projection.
\end{definition}

Many properties as in \cite{vdD} follows, for example, a
$(i_1,\ldots,i_n)$-cell has dimension $\sum_ji_j$, and if $f:X\to Y$
is a definable functions, then $\dim (X)\geq \dim (f(X))$.

\section{Preservation of balls}\label{spb}

For $o$-minimal structures, piecewise monotonicity of definable
functions plays a key role. On a general $b$-minimal structure,
there is no order $<$, so functions cannot be called monotone.
Nevertheless, the Monotonicity Theorem for $o$-minimal structures
does have an analogue for $b$-minimal structures. It is not a
consequence of $b$-minimality but has to be required as an extra
property, named preservation of (all) balls. When we look at an
$o$-minimal structure as a $b$-minimal structure as we do below,
preservation of all balls is a consequence of the Monotonicity
Theorem. The notion is especially useful for Henselian valued fields
in the context of motivic integration \cite{CLoes}, where it is used
for the change of variables in one variable, see below.

\begin{definitionm}[Preservation of balls]\label{pb}
Let $\cM$ be a $b$-minimal $\cL$-model. We say that $\cM$
\emph{preserves balls} if for every set of parameters $A$ and
$A$-definable function
 $$F:X\subset M\to M$$
there is a $A$-definable function
$$
f:X\to S
$$
as in (b1) such that for each $s\in S$
$$
F(f^{-1}(s))
$$
is either a ball or a point.

If moreover there exists such $f$ such that for every map $f_1:X\to S_1$ as in
(b1) refining $f$ (in the sense that the fibers of $f_1$ partition
the fibers of $f$) the set
 $$
F(f_1^{-1}(s_1)) $$ is also either a ball or a point for each
$s_1\in S_1$, then say that $\cM$ \emph{preserves all balls}.

We say that a $b$-minimal theory \emph{preserves balls}
(resp.~\emph{preserves all balls}) when all its models do.
\end{definitionm}

\subsection{} Let's give an example of $p$-adic integration and its
change of variables formula in one variable, using preservation of
balls.

If one integrates $|f(x)|$ over $\ZZ_p$ with $f$, say, a
semi-algebraic function $\ZZ_p\to\ZZ_p$, and $|\cdot|$ the $p$-adic
norm, it is useful to know a $b$-minimal cell decomposition of
$\ZZ_p$ relative to $\ord(f)$. That is, one takes for $X$ the
definable subset of $\ZZ_p\times (\ZZ\cup\{+\infty\})$ given by
$\ord(f(x))=a$ for $x$ in $\ZZ_p$ and $a$ in $\ZZ \cup\{+\infty\}$
and one takes a $b$-minimal cell decomposition of $X$ over
$\ZZ\cup\{+\infty\}$ to find cells $X_j$ over $\ZZ\cup\{+\infty\}$
with presentation $f_j:X_j\to (\ZZ\cup\{+\infty\})\times \ZZ^m$ for
some $m$. The fibers of $f_j$ are either balls or points, depending
on $j$ only, and since points have zero measure we can focus on
$1$-cells. Then
\begin{equation}\label{int1}
\int_{\ZZ_p}|f(x)||dx|,
\end{equation}
with $|dx|$ the Haar measure, is easily integrated, since the volume
of a ball is an easy function of its size, and since $\ord(f(x))$ by
construction is constant on the fibers of the $f_j$. Since the
measure of a ball is of the form $p^b$ for some $b\in\ZZ$, and since
$|f(x)|$ for any $x$ is of the form $p^{-b'}$ for some
$b'\in\ZZ\cup\{+\infty\}$, the integral (\ref{int1}) equals a
converging sum
\begin{equation}\label{Pint}
\sum_{a\in S} p^{-b(a)},
\end{equation}
with $S$ a Presburger set, and $b:S\to \ZZ\cup\{+\infty\}$ a
Presburger function.  Indeed, one rewrites $\ZZ_p$ as the ``tame''
disjoint union of the balls occurring  in the $1$-cells (these balls
are parameterized by a single Presburger set $S$), and on each such
ball, say parameterized by $a\in S$, one multiplies the volume of
the ball, $p^{v(a)}$, with the value of $|f(x)|=p^{-w(a)}$, where
$v$ and $w$ are Presburger functions, to obtain
$p^{-b(a)}=p^{v(a)-w(a)}$, and one then sums $p^{-b(a)}$ over $S$.

\par
In a semi-algebraic setup, preservation of balls holds such that
moreover the size of the balls is changed in a way compatible with
the Jacobian.  If $g:A\subset\ZZ_p\to \ZZ_p$ is a semi-algebraic
bijection, then
$$
\int_{\ZZ_p}|f(x)||dx|=\int_{A}|f\circ g (y)||{\rm Jac}(g)(y)||dy|,
$$
by the change of variables formula. This change of variables formula
holds here by general theory of the Haar measure on $p$-adic fields,
but such arguments fail for motivic integrals because they involve
much more general valued fields, like $k((t))$ with $k$ of
characteristic zero. However, if one takes the above cell
decomposition such that balls are preserved through $g^{-1}$ and
such that their sizes change as predicted by the Jacobian, then we
can translate both integrals into Presburger sums as (\ref{Pint})
which one sees are exactly the same Presburger sums. Indeed, the
norm of the Jacobian makes up for the difference in size of a ball
$B_a$ and its inverse image $g^{-1}(B_a)$. Thus one finds an
alternative proof of the change of variables formula in one
variable. In \cite{CLb} the motivic case is worked out. That the
preservation of balls also changes the size of the balls w.r.t. the
Jacobian, is a corollary of Weierstrass division and thus also holds
in a motivic setting and even in subanalytic motivic settings, as
long as the Henselian valued field has characteristic zero.


%

\section{Some examples of $b$-minimal structures
}\label{sectionhen}

\subsection{$o$-minimal structures and non $o$-minimal expansions}\label{somin}
Any $o$-minimal structure $R$ admits  a natural $b$-minimal
expansion by taking as main sort $R$ with the induced structure,
 the two point set $\{0,1\}$ as auxiliary sort and two constant
symbols to denote these auxiliary points. A possible interpretation
for $B$ is clear, for example,
\begin{equation*}
\begin{split}
B=\{(x,y,m)\in R^2\times R \mid x<m<y &\mbox{ when }x<y,\\ x<m
&\mbox{ when }x=y,\\\quad \mbox{ and} \ m<y &\mbox{ when } x>y
 \},
 \end{split}
\end{equation*}
so that in the $m$ variable one gets all open intervals as fibers of
$B$ above $R^2$. Property (b3) and preservation of all balls is in
this case a corollary of the Monotonicity Theorem for $o$-minimal
structures.\\

The notion of $b$-minimality leaves much more room for expansions
than the notion of $o$-minimality: some structures on the real
numbers are not $o$-minimal but are naturally $b$-minimal, for
example, the field of real numbers with a predicate for the integer
powers of $2$ are $b$-minimal by \cite{vdd2} when adding to the
above language the set of integer powers of $2$ as auxiliary sort
and the natural inclusion of it into $\RR$ as function symbol.

Recently \cite{Wilkie2}, Wilkie extended van den Dries's
construction to polynomially bounded structures, hence finding new
entire analytic functions on the reals (other than $\exp$) with tame
geometry. These structures seem to be $b$-minimal as well, w.r.t.
similar auxiliary sorts as for van den Dries's structure
$\RR,2^\ZZ$.

\subsection{Henselian valued fields of characteristic zero}
In \cite{CLb} is proved that the theory of Henselian valued fields
of characteristic zero is $b$-minimal, in a natural definitial
expansion of the valued field language, by adapting the Cohen -
Denef proof. As far as we know, this is the first written instance
of cell decomposition in mixed characteristic for unbounded
ramification.

Let $\Hens$ denote  the collection of all Henselian valued fields of
characteristic zero (hence mixed characteristic is allowed).

For $K$ in $\Hens$, write $\Ko$ for the valuation ring and $M_K$ for
the maximal ideal of $\Ko$.

For $n>0$ an integer, set $nM_K=\{nm\mid m\in M_K\}$ and consider
the natural group morphism
$$
rv_{n}:K^\times\to K^\times/1+nM_K
$$
which we extend to  $rv_n:K\to (K^\times/1+nM_K)\cup\{0\}$ by
sending $0$ to $0$.

For every $n>0$ we write $RV_{n}(K)$ for
$$RV_{n}(K):=(K^\times/1+nM_K)\cup \{0\},$$
$rv$ for $rv_1$ and $RV$ for $RV_1$.

We define the  family $B(K)$ of balls by
$$
B(K)=\{(a,b,x)\in K^\times \times K^2 \mid |x-b|<|a|\}.
$$
Hence, a ball is by definition any set of the form $B(a,b)=\{x\in
K\mid |x-b|<|a|\}$ with $a$ nonzero.

It is known that $\THens$ allows elimination of valued field
quantifiers in the language $\LHens$ by results by Scanlon
\cite{Scanlon}, F.V. Kuhlmann and Basarab \cite{KuhlBas}.

\begin{theorem}\label{cdh}
The theory $\THens$ is $b$-minimal. Moreover, $\THens$ preserves all
balls.
\end{theorem}

\begin{remark}
In fact, a slightly stronger cell decomposition theorem than Theorem
\ref{ncd} holds for $\THens$, namely a cell decomposition with
centers. We refer to \cite{CLb} to find back the full statement of
cell decomposition with centers and the definition of a center of a
cell in a $b$-minimal context.
\end{remark}

The search for an expansion of $\THens$ with a nontrivial entire
analytic function is open and challenging. Nevertheless, in
\cite{CLip} $b$-minimality for a broad class of analytic expansions
of $\THens$ is proved. This class of analytic expansions is an
axiomatization of previous work \cite{CLR}.

\section{A further study and context}

Among other things, $b$-minimality is an attempt to lay the
fundamentals of a tame geometry on Henselian valued fields that is
suitable for motivic integration, as in \cite{CLoes}. We hope to
develop this theory in future work.
One goal is to generalize the study in \cite{UdiKazh} by Hrushovski
and Kazhdan on Grothendieck
rings in a $v$-minimal context to a $b$-minimal context.\\

Theories which are $v$-minimal \cite{UdiKazh}, or $p$-minimal
\cite{Haskell} plus an extra condition, are $b$-minimal, namely, for
the $p$-minimal case, under the extra condition of existence of
definable Skolem functions. Also for $C$-minimality, some extra
conditions are needed to imply $b$-minimality. For $p$-minimality,
for example, cell decomposition lacks exactly when there are no
definable Skolem functions. A possible connection with
$d$-minimality needs to be investigated further.
\\

For notions of $x$-minimality with $x=p,C,v,o$, an expansion of a
field with an entire analytic function (other than $\exp $ on the
real field) is probably impossible, intuitively since such functions
have infinitely many zeros. In a $b$-minimal context, an infinite
discrete set does not pose any problem, see section \ref{somin} for
an example, as long as it is a ``tame" union of points. So, one can
hope for nontrivial expansions of $b$-minimal fields by entire
analytic functions, as done by Wilkie with $\exp$ and other entire
functions on the reals, see section \ref{somin}.\\

We give some open questions to end with: \\

Does a $b$-minimal $\LHens$-theory of valued fields imply that the
valued fields are Henselian?\\

As soon as the main sort $M$ is a normed field, is a definable
function $f:M^n\to M$ then automatically $C^1$, that is,
continuously
differentiable?\\

Is there a weaker condition for expansions of $\cL_B$ than
preservation of (all) balls that together with (b1), (b2) and (b3)
implies preservation of (all) balls?\\

\subsection*{Acknowledgments}Many thanks to Anand Pillay for the invitation to write this
article and to Denef, Hrushovski, and Pillay for useful advice on $b$-minimality.  \\

{\small During the writing of this paper, the author was a
postdoctoral fellow of the Fund for Scientific Research - Flanders
(Belgium) (F.W.O.) and was supported by The European Commission -
Marie Curie European Individual Fellowship with contract number HPMF
CT 2005-007121.}

\bibliographystyle{amsplain}
\bibliography{anbib}

\providecommand{\bysame}{\leavevmode\hbox to3em{\hrulefill}\thinspace}
\providecommand{\MR}{\relax\ifhmode\unskip\space\fi MR }
\providecommand{\MRhref}[2]{%
  \href{http://www.ams.org/mathscinet-getitem?mr=#1}{#2}
}
\providecommand{\href}[2]{#2}
\begin{thebibliography}{10}

\bibitem{KuhlBas}
{\c{S}}.~Basarab and F.-V. Kuhlmann, \emph{An isomorphism theorem for
  {H}enselian algebraic extensions of valued fields}, Man. Math. \textbf{77}
  (1992), no.~2-3, 113 -- 126.

\bibitem{CLip}
R.~Cluckers and L.~Lipshitz, \emph{Fields with analytic structure}, submitted,
  preprint available at http://www.dma.ens.fr/$\sim$cluckers/.

\bibitem{CLR}
R.~Cluckers, L.~Lipshitz, and Z.~Robinson, \emph{Analytic cell decomposition
  and analytic motivic integration}, Ann. Sci. \'{E}cole Norm. Sup. \textbf{39}
  (2006), no.~4, 535--568, arxiv:math.AG/0503722.

\bibitem{CLb}
R.~Cluckers and F.~Loeser, \emph{{$b$}-minimality}, J. Math. Log. \textbf{7}
  (2007), no.~2, 195--227, arXiv:math.LO/0610183.

\bibitem{CLoes}
\bysame, \emph{Constructible motivic functions and motivic integration},
  Inventiones Mathematicae \textbf{173} (2008), no.~1, 23--121,
  arxiv:math.AG/0410203.

\bibitem{Cohen}
P.~J. Cohen, \emph{Decision procedures for real and $p$-adic fields}, Comm.
  Pure Appl. Math. \textbf{22} (1969), 131--151.

\bibitem{Denef3}
J.~Denef, \emph{On the evaluation of certain $p$-adic integrals}, Th\'eorie des
  nombres, S\'emin. Delange-Pisot-Poitou 1983--84, vol.~59, 1985, pp.~25--47.

\bibitem{Denef2}
\bysame, \emph{$p$-adic semialgebraic sets and cell decomposition}, Journal
  f{\"u}r die reine und angewandte Mathematik \textbf{369} (1986), 154--166.

\bibitem{vdd2}
{L. van den} Dries, \emph{The field of reals with a predicate for the powers of
  two}, Manuscripta Math. \textbf{54} (1985), no.~1-2, 187--195.

\bibitem{vdD}
\bysame, \emph{Tame topology and o-minimal structures}, Lecture note series,
  vol. 248, Cambridge University Press, 1998.

\bibitem{HM}
D.~Haskell and D.~Macpherson, \emph{Cell decompositions of {${\rm C}$}-minimal
  structures}, Ann. Pure Appl. Logic \textbf{66} (1994), no.~2, 113--162.

\bibitem{Haskell}
\bysame, \emph{A version of o-minimality for the $p$-adics}, J. Symbolic Logic
  \textbf{62} (1997), no.~4, 1075--1092.

\bibitem{UdiKazh}
E.~Hrushovski and D.~Kazhdan, \emph{Integration in valued fields}, Algebraic
  geometry and number theory, Progr. Math., vol. 253, Birkh\"auser Boston,
  Boston, MA, 2006, pp.~261--405.

\bibitem{Miller}
C.~Miller, \emph{Tameness in expansions of the real field}, Logic Colloquium
  '01 (Urbana, IL), Lect. Notes Log., vol.~20, Assoc. Symbol. Logic, 2005,
  pp.~281--316.

\bibitem{Pas}
J.~Pas, \emph{Uniform $p$-adic cell decomposition and local zeta functions},
  Journal f\"ur die reine und angewandte Mathematik \textbf{399} (1989),
  137--172.

\bibitem{Pas2}
\bysame, \emph{Cell decomposition and local zeta functions in a tower of
  unramified extensions of a {$p$}-adic field}, Proc. London Math. Soc. (3)
  \textbf{60} (1990), no.~1, 37--67.

\bibitem{Scanlon}
T.~Scanlon, \emph{Valuation theory and its applications volume {II}}, Fields
  Institute Communications Series, ch.~Quantifier elimination for the relative
  {F}robenius, pp.~323 -- 352, {AMS}, {P}rovidence, 2003, {C}onference
  {P}roceedings of the {I}nternational {C}onference on {V}aluation {T}heory
  ({S}askatoon, 1999).

\bibitem{Wilkie2}
A.~J. Wilkie, \emph{Adding a multiplicative group to a polynomially bounded
  structure},  (October 13, 2006), Lecture at the ENS, Paris.

\end{thebibliography}

\end{document}